\documentclass[english]{amsart}

\usepackage{amsmath} 
\usepackage{amssymb} 
\usepackage{amscd}
\usepackage{tabularx}

\title{Companion forms over totally real fields, II}
\author{Toby Gee}


\newcommand{\F}{\mathbb{F}} 
\DeclareMathOperator{\Frob}{Frob} 

 \newcommand{\To}{\longrightarrow}
 \newcommand{\isoto}{\stackrel{\sim}{\To}}

 \newcommand{\bigO}{\mathcal{O}}

 \newcommand{\Z}{\mathbb{Z}}
 
 \newcommand{\Q}{\mathbb{Q}}

 \newcommand{\Gal}{\operatorname{Gal}}
 \newcommand{\GL}{\operatorname{GL}}

 \newcommand{\proj}{\operatorname{proj}}

\newcommand{\SL}{\operatorname{SL}}
\newcommand{\PSL}{\operatorname{PSL}}
\newcommand{\ad}{\operatorname{ad}}

\newcommand{\Hom}{\operatorname{Hom}}

\usepackage{amsthm}

 \newtheorem{thm}{Theorem}[section]
 
 \newtheorem{lemma}[thm]{Lemma}
 
 \theoremstyle{definition}
 
 \theoremstyle{remark}
 
\begin{document}
\maketitle
\begin{abstract}

We prove a companion forms theorem for mod $l$ Hilbert modular forms. This
work generalises results of Gross and Coleman--Voloch for modular forms over
$\Q$, and gives a new proof of their results in many cases. The methods used are completely different
to previous work in this area, and rely on modularity lifting theorems and the general theory
of deformations of Galois representations.
\end{abstract}
\section{Introduction}

If $f\in S_k(\Gamma_1(N);\overline{\F}_l)(\epsilon)$ is a mod $l$ cuspidal
eigenform, where $l\nmid N$, there is a continuous, odd, semisimple Galois
representation
\[\overline{\rho}_f:\Gal(\overline{\Q}/\Q)\To\GL_2(\overline{\F}_l)
\]
attached to $f$. A famous conjecture of Serre predicts that all continuous
odd irreducible mod $l$ representations should arise in this fashion.
Furthermore, the ``strong Serre conjecture'' predicts a minimal weight
$k_{\overline{\rho}}$ and level $N_{\overline{\rho}}$, in the sense that $\overline{\rho} \cong\overline{\rho}_g$ for
some eigenform $g$ of weight $k_{\overline{\rho}}$ and level $N_{\overline{\rho}}$ (prime to $l$),
and if $\overline{\rho}\cong\overline{\rho}_f$ for some eigenform $f$ of weight $k$ and level $N$
prime to $l$ then $N_{\overline{\rho}}|N$ and $k\geq k_{\overline{\rho}}$. The question as to
whether all continuous odd irreducible mod $l$ Galois representations are
modular in this sense is still open, although it is expected to be known soon as a result of the work of Khare and
Wintenberger (see for example \cite{kw06}). However the implication ``weak Serre
$\Rightarrow$ strong Serre'' was essentially established over ten years ago (aside from a few cases
where $l=2$).

In solving the problem of weight optimisation it becomes necessary to
consider the companion forms problem; that is, the question of when it can
occur that we have $f=\sum a_nq^n$ of weight $2\leq k\leq l$ with $a_l\neq
0$, and an eigenform $g=\sum b_n q^n$ of weight $k'=l+1-k$ such that
$na_n=n^k b_n$ for all $n$. Serre conjectured that this can occur if and
only if the representation $\rho_f$ is tamely ramified above $l$. This
conjecture has been settled in most cases in the papers of Gross
(\cite{gro90}) and Coleman-Voloch (\cite{cv92}).

Our earlier paper \cite{gee04} generalised these results to the case of
parallel weight Hilbert modular forms over totally real fields $F$ in which
$l$ splits completely, by generalising the methods of \cite{cv92}. In this
paper we take a completely different and rather more conceptual approach; we
construct our companion form by using a method of Ramakrishna to find an
appropriate characteristic zero Galois representation, and then use recent
work of Kisin (\cite{kis04}) to prove that the representation is modular.
Note that our companion form is not necessarily of minimal prime-to-$l$
level, but that this is irrelevant for applications to Artin's conjecture,
and that in many cases a form of minimal level may be obtained from ours by
the methods of \cite{jar99}, \cite{sw01}, \cite{raj01} and \cite{fuj99}. In
the case of weight $l$ forms, we avoid potential difficulties with weight
$1$ forms by constructing a companion form in weight $l$.

\section{Statement of the main results}
Let $l>2$ be a prime, and let $F$ be a totally real
field. We assume that if $l>3$, $[F(\zeta_l):F]>2$ (note that this is
automatic if $l$ is unramified in $F$). Let $\epsilon$ denote both the
$l$-adic and mod $l$ cyclotomic characters; this should cause no confusion.
Let $\rho:G_K\to\GL_2(\bigO)$ be a continuous representation, where is $K$ a
finite extension of $\Q_l$, and $\bigO$ is the ring of integers in a finite
extension of $\Q_l$. We say that $\rho$ is \emph{ordinary} if it is
Barsotti-Tate, coming from an $l$-divisible group which is an extension of
an \'{e}tale group by a multiplicative group, each of rank one as
$\bigO$-modules. We say that it is \emph{potentially ordinary} if it becomes
ordinary upon restriction to an open subgroup of $G_K$. We say that a
Hilbert modular form of parallel weight 2 is \emph{(potentially) ordinary}
at a place $v|l$ if its associated Galois representation is (potentially)
ordinary at $v$. These definitions agree with those in \cite{kis04}; they
are slightly non-standard, but if the level is prime to $l$ then
this is equivalent to the $U_v$-eigenvalue being an $l$-adic unit. We say
that a Hilbert modular form of parallel weight $k$, $3\leq k\leq l$, is
\emph{ordinary} at a place $v|l$ if its $U_v$-eigenvalue is an $l$-adic
unit. Finally, we say that a modular form is \emph{(potentially) ordinary}
if it is (potentially) ordinary at all places $v|l$.

Our main theorem is the following:

\begin{thm}\label{maink}Let $g$ be an ordinary Hilbert modular eigenform of parallel weight $k$,
$2\leq k\leq l$, and level coprime to $l$. Let its associated Galois
representation be $\rho_g:G_F\to\GL_2(\overline{\Q}_l)$, so that (by
Theorem 2 of \cite{wil90}) we have, for all places $v|l$, $$\rho_g|_{G_v}\simeq\left(%
\begin{array}{cc}
  \epsilon^{k-1}\psi_{v,1} & * \\
  0 & \psi_{v,2} \\
\end{array}%
\right)$$for unramified characters $\psi_{v,1}$, $\psi_{v,2}$. Suppose that
the residual representation $\overline{\rho}_g:G_F\to\GL_2(\overline{\F}_l)$
is absolutely irreducible. Assume further that for all $v|l$  we have that
$\epsilon^{k-1}\overline{\psi}_{v,1}\neq\overline{\psi}_{v,2}$, and that the
representation
$\overline{\rho}_g|_{G_v}$ is tamely ramified, so that $$\overline{\rho}_g|_{G_v}\simeq\left(%
\begin{array}{cc}
  \epsilon^{k-1}\overline{\psi}_{v,1} & 0 \\
  0 & \overline{\psi}_{v,2} \\
\end{array}%
\right).$$ Assume in addition that if $\epsilon^{k-2}\overline{\psi}_{v,1}=\overline{\psi}_{v,2}$, then
the absolute ramification index of $F_v$ is less than $l-1$. If $k=l$ then let $k'=l$, and otherwise let $k'=l+1-k$. Then
there is a Hilbert modular form $g'$ of parallel weight $k'$ and level
coprime to $l$ satisfying
$$\overline{\rho}_{g'}\simeq\overline{\rho}_g\otimes\epsilon^{k'-1}$$and
the $U_v$-eigenvalue of $g'$ is a lift of $\overline{\psi}_{v,1}(\Frob_v)$.
\end{thm}

In fact, we work throughout with forms of parallel weight $2$, and we use
Hida theory to treat forms of more general (parallel) weight. In the case
where $\overline{\rho}_g(G_F)$ is soluble the Langlands-Tunnell theorem
makes the proof straightforward, so we concentrate on the insoluble case,
where we prove:

\begin{thm}\label{wt2}Let $\overline{\rho}_f:G_F\to\GL_2(\overline{\F}_l)$ be
an absolutely irreducible modular representation, coming from a Hilbert
eigenform $f$ of parallel weight $2$, with associated Galois representation
$\rho_f:G_F\to\GL_2(\overline{\Q}_l)$. Suppose that $\overline{\rho}_f(G_F)$
is insoluble. Suppose also that for every place $v$ of $F$ dividing $l$,
$\rho_f|_{G_{v}}$ is potentially ordinary, and
we have $$\overline{\rho}_f|_{G_v}\simeq\left(%
\begin{array}{cc}
  \epsilon^{k-1}\overline{\psi}_{v,1} & 0 \\
  0 & \overline{\psi}_{v,2} \\
\end{array}%
\right)$$where $\overline{\psi}_{v,1}$, $\overline{\psi}_{v,2}$ are
unramified characters, with $\epsilon^{k-1}\overline{\psi}_{v,1} \neq
\overline{\psi}_{v,2}$. Assume in addition that if $\epsilon^{k-2}\overline{\psi}_{v,1}=\overline{\psi}_{v,2}$, then
the absolute ramification index of $F_v$ is less than $l-1$.

If $k=l$ then let $k'=l$, and otherwise let $k'=l+1-k$. Then there is an
eigenform $f'$ of parallel weight $2$ which is potentially ordinary at all
places $v|l$ such that the mod $l$ Galois representation
$\overline{\rho}_{f'}$ associated to $f'$ satisfies
$$\overline{\rho}_{f'}\simeq\overline{\rho}_f\otimes\epsilon^{k'-1},$$
and such that at all places $v|l$ we have
$$\rho_{f'}|_{G_v}\simeq\left(%
\begin{array}{cc}
  \epsilon\omega^{k'-2}\psi_{v,2} & * \\
  0 & \psi_{v,1} \\
\end{array}%
\right)$$ with $\psi_{v,i}$ an unramified lift of $\overline{\psi}_{v,i}$
for $i$=1, 2, and $\omega$ the Teichmuller lift of $\epsilon$.
\end{thm}
\section{Lifting theorems}
Firstly, we prove a straightforward generalisation of the results of
\cite{ram02} and \cite{tay03} to totally real fields. We fix the determinant of our deformations throughout this section
(that is, we work with $\ad^0\overline{\rho}$ rather than $\ad\overline{\rho}$). We originally did this to follow \cite{tay03}, and we would like to thank the anonymous referee for
pointing out to us that, unlike over $\Q$, it is in fact necessary to do this when working over totally real fields, as
there may no longer be a choice of global determinant satisfying given local conditions. We begin by analysing
the local representation theory at primes not dividing $l$. The next lemma
is essentially contained in \cite{dia972}:

\begin{lemma}\label{localnotl}Let $p\neq l$ be a prime, and let $K$ be a finite
extension of $\Q_p$. Let $I_K$ denote the inertia subgroup of $G_K$. Let
$\sigma:G_K\to\GL_2(\F)$ be a continuous representation, with $\F$ a finite
field of characteristic $l$, and assume that $l|\#\sigma(I_K)$.

Then either $p=2$, $l=3$, and $\proj\sigma(G_K)\simeq A_4$ or $S_4$,
or $$\sigma\simeq\left(%
\begin{array}{cc}
  \epsilon\overline{\chi} & * \\
  0 & \overline{\chi} \\
\end{array}%
\right)$$with respect to some basis for some character $\overline{\chi}$.
\end{lemma}
\begin{proof}Note that $l|\#\sigma(I_K)$ if and only if
$l|\#\proj\sigma(I_K)$. We must have $\sigma|_{I_K}$ indecomposable. If
$\sigma$ is reducible, then $\sigma$ is a twist of a representation
$\bigl(\begin{smallmatrix}\psi&u\\0&1\end{smallmatrix}\bigr)$ for some
character $\psi$, with $u$ a cocycle representing a class in
$\operatorname{H}^1(G_K,k(\psi))$ whose image in
$\operatorname{H}^1(I_K,k(\psi))^{G_K}$ is non-zero; but the latter group is
zero unless $\psi=\epsilon$.

If instead $\sigma$ is irreducible but $\sigma|_{I_K}$ is reducible, then
$\sigma|_{I_K}$, being indecomposable, must fix precisely one element of
$\mathbb{P}^1(k)$. But then $\sigma$ would also have to fix this element, a
contradiction.

Assume now that $\sigma|_{I_K}$ is irreducible, and that $\sigma|_{P_K}$ is
reducible, where $P_K$ is the wild inertia subgroup of $I_K$. Then $P_K$
must fix precisely two elements of $\mathbb{P}^1(k)$ (as $\sigma|_{I_K}$ is
irreducible), so $\sigma$ is induced from a character on a ramified
quadratic extension of $K$, and thus $\sigma(I_K)$ has order $2p^r$ for some
$r\geq1$, a contradiction.

Finally, if $\sigma|_{P_K}$ is irreducible we must have $p=2$. That
$\proj\sigma(G_K)\simeq A_4$ or $S_4$ follows from the same argument as in
the proof of Proposition 2.4 of \cite{dia972}. That $l=3$ follows from
$l|\#\sigma(I_K)$.
\end{proof}

Let $\overline{\rho}:G_F\to\GL_2(\F)$ be continuous, odd, and absolutely
irreducible, with $k$ a finite field of characteristic $l$. Let $S$ denote a
finite set of finite places of $F$ which contains all places dividing $l$
and all places where $\overline{\rho}$ is ramified, and let $G_S$ denote the
Galois group of the maximal extension of $F$ unramified outside $S$. A
\emph{deformation} of $\overline{\rho}$ is a complete noetherian local ring
$(R,\mathfrak{m})$ with residue field $k$ and a continuous representation
$\rho:G_S\to\GL_2(R)$ such that $(\rho$ mod $\mathfrak{m})=\overline{\rho}$
and $\epsilon^{-1}\det\rho$ has finite order prime to $l$. We define
deformations of $\overline{\rho}|_{G_v}$ in a similar fashion.

Suppose that for each $v\in S$ we have a pair $(\mathcal{C}_v,L_v)$
satisfying the properties P1-P7 listed in section 1 of \cite{tay03}. Define
$\operatorname{H}^1_{\{L_v\}}(G_S,\ad^0\overline{\rho})$ and
$\operatorname{H}^1_{\{L^\bot_v\}}(G_S,\ad^0\overline{\rho})$ in the usual
way.

\begin{lemma}\label{32}If
$\operatorname{H}^1_{\{L^\bot_v\}}(G_S,\ad^0\overline{\rho}(1))=(0)$ then there
is an $S$-deformation $(W(k),\rho)$ of $\overline{\rho}$ such that for all
$v\in S$ we have $(W(k),\overline{\rho}|_{G_v})\in\mathcal{C}_v$.
\end{lemma}
\begin{proof}Identical to the proof of Lemma 1.1 of \cite{tay03}.
\end{proof}
\begin{lemma}\label{lift}Suppose that $\sum_{v\in S}\dim L_v\geq\sum_{v\in
S\cup\{\infty\}}\dim\operatorname{H}^0(G_v,\ad^0\overline{\rho})$. Then we
can find a finite set of places $T\supset S$ and data $(\mathcal{C}_v,L_v)$
for $v\in T-S$ satisfying conditions P1-P7 and such that
$\operatorname{H}^1_{\{L^\bot_v\}}(G_T,\ad^0\overline{\rho}(1))=(0)$.
\end{lemma}
\begin{proof}The proof of this lemma is almost identical to that of
Lemma 1.2 of \cite{tay03}. We sketch a few of the less obvious details. In
the case $l=5$, $\ad^0\overline{\rho}(G_F)\simeq A_5$, we choose $w\notin S$
such that $\mathbb{N}w\equiv 1 \text{ mod }5$ and
$\ad^0\overline{\rho}(\Frob_w)$ has order 5 (such a $w$ exists by
Cebotarev's theorem). Adding $w$ to $S$ with the pair $(\mathcal{C}_w,L_w)$
of type E3 (see below), we may assume
$\operatorname{H}^1_{\{L^\bot_v\}}(G_S,\ad^0\overline{\rho})\cap\operatorname{H}^1(\ad^0\overline{\rho}(G_F),\ad^0\overline{\rho})=(0)$.

From here on, almost exactly the same argument as in \cite{tay03} applies,
the only difference being that one must replace every occurence of ``$\Q$''
with ``$F$''. Let $M=F(\ad^0\overline{\rho},\mu_l)$. The argument is
essentially formal once one knows that there is an element
$\sigma\in\Gal(M/F)$ such that $\ad^0\overline{\rho}(\sigma)$ has an
eigenvalue $\epsilon(\sigma)\not\equiv 1\text{ mod }l$, that
$\ad^0\overline{\rho}$ is absolutely irreducible, and that
$\ad^0\overline{\rho}$ is not isomorphic to $(\ad^0\overline{\rho})(1)$. All
of these assertions follow from our assumption that $[F(\zeta_l):F]>2$ if
$l>3$, with the proofs being similar to those in \cite{ram99} (note that one
may replace the assumption that $\overline{\rho}(G_\Q)\supseteq\SL_2(\F)$ in
\cite{ram99} with the assumption that
$\proj\overline{\rho}(G_\Q)\supseteq\PSL_2(\F)$ without affecting the
proofs). For example, to check that $\ad^0\overline{\rho}$ is not isomorphic
to $(\ad^0\overline{\rho})(1)$ it is enough to prove that there is an
element $\sigma'\in\Gal(M/F)$ such that all of the eigenvalues of
$\ad^0\overline{\rho}(\sigma')$ are $1$, and $\epsilon(\sigma')\neq 1$. The existence
of $\sigma$ and $\sigma'$ follows exactly as in the proof of Theorem 2 of
\cite{ram99}.\end{proof}

We now give examples of pairs $(\mathcal{C}_v,L_v)$. Again, our pairs are
very similar to those in section 1 of \cite{tay03}, and the verification of
the required properties is almost identical. We use the notation of
\cite{tay03} for ease of comparison with that paper.

\begin{itemize}\item E1. Suppose that $v\nmid l$ and that
$l\nmid\#\overline{\rho}(I_v)$. Take $\mathcal{C}_v$ to be the class of
lifts of $\overline{\rho}|_{G_v}$ which factor through
$G_v/(I_v\cap\ker\overline{\rho})$ and let $L_v$ be
$\operatorname{H}^1(G_v/I_v,(\ad^0\overline{\rho})^{I_v})$. Then it is
straightforward to see that properties P1-P7 are satisfied, and
that\begin{itemize}\item
$\operatorname{H}^2(G_v/(I_v\cap\ker\overline{\rho}),\ad^0\overline{\rho})\simeq\operatorname{H}^2(G_v/I_v,(\ad^0\overline{\rho})^{I_v})=(0)$,
(as $G_v/I_v\simeq\hat{\Z}$ has cohomological dimension 1),
\item
$\operatorname{H}^1(G_v/(I_v\cap\ker\overline{\rho}),\ad^0\overline{\rho})=L_v\subset
\operatorname{H}^1(G_v,\ad^0\overline{\rho})$, \item $\dim L_v=\dim
\operatorname{H}^0(G_v,\ad^0\overline{\rho})$ (by the local Euler
characteristic formula).\end{itemize}

\item E2. (Note that our definitions here differ slightly from those in \cite{tay03}; we thank Richard Taylor for
explaining this modification to us.) Suppose that $l=3$, that $v|2$, and
that $(\ad^0\overline{\rho})(G_v)\isoto S_4$. Take $\mathcal{C}_v$ to be the
class of lifts of $\overline{\rho}|_{G_v}$ which factor through
$G_v/(I_v\cap\ker\overline{\rho})$ and let $L_v$ be
$\operatorname{H}^1(G_v/I_v,(\ad^0\overline{\rho})^{I_v})$. The verification
of properties P1-P7 is then as in \cite{tay03}, except that to check that
$\operatorname{H}^i(\overline{\rho}(I_v),\ad^0\overline{\rho})=(0)$ for all
$i\geq 0$ one uses the Hochschild-Serre spectral sequence and the fact that
$\operatorname{H}^i(C_2 \times C_2,\ad^0 \overline{\rho})=(0)$ for all $i\geq 0$.

\item E3. Suppose that $v\nmid l$, that either $\mathbb{N}v\not\equiv 1$ (mod $l$) or $l|\#\overline{\rho}(G_v)$, and that with respect to some basis $e_1$, $e_2$ of
$\F^2$ the restriction $\overline{\rho}|_{G_{v}}$ has the form

$$\left(%
\begin{array}{cc}
  \epsilon\overline{\chi} & * \\
  0 & \overline{\chi} \\
\end{array}%
\right).$$Take $\mathcal{C}_v$ to be the class of deformations of the form
(with respect to some basis)
$$\left(%
\begin{array}{cc}
  \epsilon\chi & * \\
  0 & \chi \\
\end{array}%
\right)$$with $\chi$ lifting $\overline{\chi}$, and take $L_v$ to be the
image of
$$\operatorname{H}^1(G_v,\Hom(\F e_2,\F e_1))\to\operatorname{H}^1(G_v,\ad^0\overline{\rho}).$$

That the pair $(\mathcal{C}_v,L_v)$ satisfies the properties P1-P7 follows
from an identical argument to that in \cite{tay03}. An identical calculation
to that in \cite{tay03} shows that $\dim
L_v=\dim\operatorname{H}^0(G_v,\ad^0\overline{\rho})$.

\item E4. Suppose that $v|l$ and that with respect to some basis
$e_1$, $e_2$ of $\F^2$, $\overline{\rho}|_{G_v}$ has the form $$\left(%
\begin{array}{cc}
  \epsilon\overline{\chi_1} & 0 \\
  0 & \overline{\chi_2} \\
\end{array}%
\right).$$ Suppose also that $\overline{\chi}_1\neq\overline{\chi}_2$ and
that $\epsilon\overline{\chi}_1\neq\overline{\chi}_2$. Take
$\mathcal{C}_v$ to consist of all deformations of the form $$\left(%
\begin{array}{cc}
  \epsilon\chi_1 & * \\
  0 & \chi_2 \\
\end{array}%
\right)$$where $\chi_1$, $\chi_2$ are tamely ramified lifts of
$\overline{\chi}_1$, $\overline{\chi}_2$ respectively. Let
$U^0=\Hom(\F e_2,\F e_1)$, and let $L_v$ be the kernel of the map
$\operatorname{H}^1(G_v,\ad^0\overline{\rho})\to\operatorname{H}^1(I_v,\ad^0\overline{\rho}/U^0)^{G_v/I_v}$.
The verification of properties P1-P7 follows as in \cite{tay03}, and we may
compute $\dim L_v$ via a similar computation to that in the proof of Lemma 5
of \cite{ram02}.

Note firstly that by local duality and the assumption that
$\overline{\chi}_1\neq\overline{\chi}_2$ we have
$\operatorname{H}^2(G_v,U^0)=0$. Thus the short exact sequence
$$0\to U^0\to\ad^0\overline{\rho}\to\ad^0\overline{\rho}/U^0\to 0$$yields
an exact sequence
$$\operatorname{H}^1(G_v,\ad^0\overline{\rho})\to\operatorname{H}^1(G_v,\ad^0\overline{\rho}/U^0)\to
0.$$ Inflation-restriction gives us an exact sequence
$$0\to\operatorname{H}^1(G_v/I_v,(\ad^0\overline{\rho}/U^0)^{I_v})\to\operatorname{H}^1(G_v,\ad^0\overline{\rho}/U^0)\to\operatorname{H}^1(I_v,\ad^0\overline{\rho}/U^0)^{G_v/I_v}\to
0,$$ and combining these two sequences shows that the map
$\operatorname{H}^1(G_v,\ad^0\overline{\rho})\to\operatorname{H}^1(I_v,\ad^0\overline{\rho}/U^0)^{G_v/I_v}$
is surjective. Thus \begin{align*}\dim L_v
&=\dim\operatorname{H}^1(G_v,\ad^0\overline{\rho})-
\dim\operatorname{H}^1(I_v,\ad^0\overline{\rho}/U^0)^{G_v/I_v}\\
&=\dim\operatorname{H}^1(G_v,\ad^0\overline{\rho})-\dim\operatorname{H}^1(G_v,\ad^0\overline{\rho}/U^0)+\dim\operatorname{H}^1(G_v/I_v,(\ad^0\overline{\rho}/U^0)^{I_v})\\
&=\dim\operatorname{H}^1(G_v,\ad^0\overline{\rho})-\dim\operatorname{H}^1(G_v,\ad^0\overline{\rho}/U^0)\\&\
+\dim\operatorname{H}^0(G_v,\ad^0\overline{\rho}/U^0)\
\text{(by Lemma 3 of \cite{ram02})}\\
&=\dim\operatorname{H}^0(G_v,\ad^0\overline{\rho})+\dim\operatorname{H}^2(G_v,\ad^0\overline{\rho})-\dim\operatorname{H}^2(G_v,\ad^0\overline{\rho}/U^0)\\&\
+[F_v:\Q_l]\ \text{(local Euler
characteristic)}\\&=[F_v:\Q_l]+\dim\operatorname{H}^0(G_v,\ad^0\overline{\rho}).\end{align*}

\item{BT.} Suppose that $v|l$ and that with respect to some basis
$e_1$, $e_2$ of $\F^2$, $\overline{\rho}|_{G_v}$ has the form $$\left(%
\begin{array}{cc}
  \epsilon\overline{\chi} & 0 \\
  0 & \overline{\chi} \\
\end{array}%
\right)$$for some unramified character $\overline{\chi}$. Assume also that
$\epsilon$ is not trivial (that is, that $F_v$ does not contain $\Q_l(\zeta_l)$).
Take
$\mathcal{C}_v$ to consist of all flat deformations of the form $$\left(%
\begin{array}{cc}
  \epsilon\chi_1 & * \\
  0 & \chi_2 \\
\end{array}%
\right)$$where $\chi_1$, $\chi_2$ are unramified lifts of $\overline{\chi}$,
 Then it follows from Corollary 2.5.16 of \cite{kis04} that
there is an $L_v$ of dimension
$[F_v:\Q_l]+\dim\operatorname{H}^0(G_v,\ad^0\overline{\rho})$ so that
properties P1-P7 are all satisfied.

\end{itemize}
Set $\overline{\rho}=\overline{\rho}_f\otimes\epsilon^{k'-1}$. We are now in
a position to prove:
\begin{thm}\label{btlift}There is a deformation $\rho$ of $\overline{\rho}$ to $W(k)$
such that at all places $v|l$ we have $\rho|_{G_v}$ potentially ordinary,
and
$$\rho|_{G_v}\simeq\left(%
\begin{array}{cc}
  \epsilon\omega^{k'-2}\psi_{v,2} & * \\
  0 & \psi_{v,1} \\
\end{array}%
\right)$$ with $\psi_{v,i}$ an unramified lift of $\overline{\psi}_{v,i}$
for $i$=1, 2, and $\omega$ the Teichm\"{u}ller lift of $\epsilon$.
\end{thm}

\begin{proof}This follows almost at once from Lemmas \ref{32} and \ref{lift}. By
Lemma \ref{localnotl} we can choose $(\mathcal{C}_v,L_v)$ for all $v\nmid
l$, with $\dim L_v=\dim \operatorname{H}^0(G_v,\ad^0\overline{\rho})$
(simply choose as in examples E1 or E3). At places $v|l$, we choose
$(\mathcal{C}_v,L_v)$ as in examples E4 or BT, so that $\dim
L_v=[F_v:\Q_l]+\dim\operatorname{H}^0(G_v,\ad^0\overline{\rho})$. Then as
$\sum_{v|l}[F_v:\Q_l]=[F:\Q]$, we have $\sum_{v\in S}\dim L_v=\sum_{v\in
S\cup\{v|\infty\}}\dim\operatorname{H}^0(G_v,\ad^0\overline{\rho})$, so a
deformation as in Lemma \ref{32} exists. That the $\psi_{v,i}$ are
unramified follows from the fact that they are tamely ramified lifts of
unramified characters.

It remains to check that $\rho|_{G_v}$ is potentially ordinary.  By the
remarks in section 2.4.15 of \cite{kis04} it suffices to
check that it is potentially Barsotti-Tate. This is immediate if we are in
the case BT, so suppose we are considering deformations as in E4. By the
proposition in section 3.1 of \cite{per94}, $\rho|_{G_v}$ is potentially
semistable, and it clearly has Hodge-Tate weights in $\{0,1\}$, so by
Theorem 5.3.2 of \cite{bre00} it suffices to check that it is potentially
crystalline. In order to check this, we consider the Weil-Deligne
representation $WD(\rho|_{G_v})$ (see Appendix B of \cite{cdt} for the
definition of $WD(\sigma)$ for any potentially semistable $p$-adic
representation $\sigma$ of $G_v$). We need to check that the associated
nilpotent endomorphism $N$ is zero. As is well-known, $N=0$ unless
$WD(\rho|_{G_v})$ is a twist of the Steinberg representation, which cannot
happen because of our assumption that we are not in the BT case.
\end{proof}

Theorem \ref{wt2} now follows immediately from:

\begin{thm}The representation $\rho$ is modular.
\end{thm}

\begin{proof}This is an easy application of Theorem 3.5.5 of
\cite{kis04}. We need to check that $\rho$ is strongly residually modular. The representation $\rho_f\otimes\omega^{k'-1}$ (where $\omega$ is the
Teichm\"{u}ller lift of $\epsilon$) is certainly modular, with residual
representation $\overline{\rho}$, and $\overline{\rho}(G_F)$ is insolvable. Furthermore, it is automatically potentially ordinary at all places $v|l$ with $\epsilon^{k-2}\overline{\psi}_{v,1}\neq\overline{\psi}_{v,2}$.
By Theorem 6.2 of \cite{jar04} and our assumption that if $\epsilon^{k-2}\overline{\psi}_{v,1}=\overline{\psi}_{v,2}$ the absolute ramification
index of $F_v$ is less than $l-1$, we may replace $\rho_f\otimes\omega^{k'-1}$ with a modular lift of $\overline{\rho}$
which is potentially ordinary at all places $v|l$. By construction, $\rho$ is potentially ordinary at all places $v|l$, so we are done.
\end{proof}

We now prove Theorem \ref{maink}.
Firstly, suppose that $\overline{\rho}_g(G_F)$ is insoluble. Then Hida
theory (see \cite{wil90} or \cite{hid88}) provides us with a
weight $2$ form $f$ which satisfies the hypotheses of Theorem \ref{wt2}, and
which has $\overline{\rho}_f\simeq\overline{\rho}_g$ (that $f$ is
potentially ordinary follows as in the proof of Theorem \ref{btlift}). Then
Theorem \ref{wt2} provides us with a Hilbert modular form $f'$ of parallel
weight $2$ with
$\overline{\rho}_{f'}\simeq\overline{\rho}_f\otimes\epsilon^{k'-1}$
and $$\rho_{f'}|_{G_v}\simeq\left(%
\begin{array}{cc}
  \epsilon\omega^{k'-2}\psi_{v,2} & * \\
  0 & \psi_{v,1} \\
\end{array}%
\right)$$for all places $v|l$, with $\psi_{v,1}$ an unramified lift of
$\overline{\psi}_{v,1}$. Then Lemma 3.4.2 of \cite{kis04} shows that $f'$ has $U_v$-eigenvalue
$\psi_{v,1}(\Frob_v)$, an $l$-adic unit. The existence of $g'$ now follows
from Hida theory.

Now suppose that $\overline{\rho}_f(G_F)$ is solvable. Then there is a lift
of $\overline{\rho}_f\otimes\epsilon^{k'-1}$ to a characteristic zero
representation, which comes from a Hilbert modular form of parallel weight 1
by the Langlands-Tunnell theorem. Such a form is necessarily ordinary in the
sense of Hida theory, and the theorem follows by Hida theory as in the
insoluble case.

\def\cprime{$'$}
\providecommand{\bysame}{\leavevmode\hbox to3em{\hrulefill}\thinspace}
\providecommand{\MR}{\relax\ifhmode\unskip\space\fi MR }
\providecommand{\MRhref}[2]{%
  \href{http://www.ams.org/mathscinet-getitem?mr=#1}{#2}
}
\providecommand{\href}[2]{#2}

\end{document}